# 3D Adaptive Central Schemes: part I Algorithms for Assembling the Dual Mesh

S. Noelle,[†] W. Rosenbaum,[†] and M. Rumpf[‡]

October 1, 2018
supported by DFG-SPP, ANumE

**Abstract.** Central schemes are frequently used for incompressible and compressible flow calculations. The present paper is the first in a forthcoming series where a new approach to a 2nd order accurate Finite Volume scheme operating on cartesian grids is discussed. Here we start with an adaptively refined cartesian *primal grid* in 3D and present a construction technique for the staggered *dual grid* based on $L^\infty$-Voronoi cells. The local refinement constellation on the primal grid leads to a finite number of uniquely defined *local patterns* on a primal cell. Assembling adjacent local patterns forms the dual grid. All local patterns can be analysed in advance. Later, running the numerical scheme on staggered grids, all necessary geometric information can instantly be retrieved from lookup-tables. The new scheme is compared to established ones in terms of algorithmical complexity and computational effort.

## 1 Introduction

CFD simulations, especially in 3D, produce large amounts of data. In order to minimize memory requirements and computing time without sacrificing high spatial resolution in regions of interest modern numerical schemes are usually based on adaptive grids. Because of their simple structure and the ease of data access, cartesian grids are particularly popular.

Staggered grids are a pair of meshes of the same computational domain whose nodes (in 1D), edges (in 2D) and faces (in 3D) do not coincide. While the shape of a staggered grid is canonical on an uniform mesh it becomes rather complicated for an underlying adaptively refined grid, in particular in three space dimensions. In the present paper we suggest a construction technique for a staggered grid which is dual to an adaptive cartesian grid.

Our focus of interest is the class of central schemes for hyperbolic conservation laws, which became quite popular during the last decade. In this case, staggered grids circumvent the possibly costly numerical solution of Riemann problems. The prototype of all central schemes is the first-order Lax-Friedrich scheme [?]. Higher order extensions were first introduced by Nessyahu and Tadmor [?], and developed further by many others (see [?, ?, ?, ?, ?, ?] and the references therein).

---

[†]RWTH Aachen, Germany
[‡]University of Duisburg-Essen, Germany

Let us briefly summarize our staggered grid construction. The cells of the dual grid are defined as the $L^\infty$-Voronoi regions around the vertices of the primal grid. Even though it is in principle simple to construct these Voronoi regions, it becomes a very complex task if the primal grid is adaptively refined, especially in three dimensions. Our approach, motivated by previous work on the 2D case [?], is based on a local decomposition of each primal cell into dual subcells. We call this decomposition of a primal cell the *local pattern*. The dual subcells are later assembled at the vertices. Since we use a cartesian grid with graded refinement, the number of different local decompositions is bounded. In [?], we used Polya theory to show that there are 227 combinatorially different local patterns in 3D. All these patterns can be analysed and stored in advance. Later, running the numerical scheme on staggered grids, all necessary geometric information can instantly be retrieved from lookup-tables.

A popular alternative dual grid construction on adaptive cartesian grids is based on "diamond cells", developed for example in [?, ?]. Although this geometrical description is very simple, we show in Section 2 that diamond grids increase the numerical cost of the finite volume scheme considerably.

The paper is organized as follows: In Section 2 we investigate Voronoi decompositions with respect to several norms and compare them with diamond dual grids. In order to explain the techniques as clearly as possible, we begin by constructing the local patterns for 2D grids in Section 3. In Section 4, we generalize the pattern construction to 3D grids. Section 5 presents first numerical applications in 3D, discusses algorithmical complexity and poses some open questions to be treated in forthcoming papers.

## 2 Problem setting and concept

In this section we introduce the necessary notation, formulate properties of the primal and the corresponding dual grid and explain the basic idea of the dual grid construction. We discuss both the popular "diamond cell construction" as well as the Voronoi decomposition with respect to different norms in some detail and argue for our favorite choice.

Throughout this paper we use the capital letters $G$, $C$, $F$, $E$ and $N$ for a *grid, cell, face, edge* resp. *node*, capital calligraphic letters to denote *sets*, like $\mathcal{N} = \bigcup N$, subscripts $p$ and $d$ to distinguish between *primal* and *dual* objects, an asterisk $^*$ to denote objects from the corresponding staggered grid (the dual grid is staggered to the primal grid and the converse), subscripts $i$, $j$ and $k$ for node indices, and the superscript $^+$ to label local objects on single cells.

Our primary interest of constructing dual grids is directed towards the numerical solution of systems of conservation laws

$$(1) \qquad u_t(x,t) + \operatorname{div} F(u(x,t)) = 0$$

in three space dimensions. Here, $x \in \mathbf{R}^3$ denotes the space variable, $t$ the evolution time, $u(x,t) \in \mathbf{R}^m$ the solution vector in terms of $m$ conservative variables, and $F = (f_i)_{i=1}^3$ the vector of the three directional flux functions $f_i : \mathbf{R}^m \to \mathbf{R}^m$, $i = 1, 2, 3$.

The final goal is to implement a 3D-extension of the second order accurate central scheme proposed in [?]. A dimensionless formulation of the finite volume scheme reads as follows:

Let $v = v(x,t)$ be a cellwise smooth approximate solution of (1) with cell averages

$$(2) \qquad \bar{v}_C^n = \frac{1}{|C|} \int_C v(x, t^n)\, dx.$$

The standard finite volume update is

$$(3) \qquad \bar{v}_C^{n+1} = \bar{v}_C^n + \frac{1}{|C|} \int_{t^n}^{t^{n+1}} \int_{\partial C} F(v) \cdot \mathfrak{n} \, dx\, dt$$

Note that the flux integral should be replaced by a quadrature rule, and the approximate solution $v$ needs to be extrapolated in time (for a fully discrete scheme). Since the solution $v$ may be discontinuous at the cell boundary $\partial C$, one would have to replace the flux function $F(v)$ by a Riemann solver. This can be avoided when using staggered grids: now mean values on the timelevel $t^{n+1}$ are computed on cells $C^*$ of the corresponding staggered grid via

$$(4) \qquad \bar{v}_{C^*}^{n+1} = \frac{1}{|C^*|} \sum_{C \cap C^* \neq \emptyset} \left\{ \underbrace{\int_{C \cap C^*} v(x, t^n)\, dx}_{I_1} + \underbrace{\int_{t^n}^{t^{n+1}} \int_{C \cap \partial C^*} F(v) \cdot \mathfrak{n}\, dx\, dt}_{I_2} \right\}$$

As before, one needs quadrature rules for the integrals and time extrapolation for $v$. For the flux integral $I_2$ this is possible in such a way that $F(v)$ needs to be computed only at points $(x,t)$ where $v$ is continuous [?, ?, ?]. Therefore, there is no need to solve Riemann problems for central schemes. Still, the evaluation of the flux function may be the most expensive part of the finite volume scheme. Therefore, the numerical scheme should minimize the number of quadrature points in $I_2$.

The cartesian *primal grid* $G_p$ is supposed to be an adaptively refined 3D cartesian mesh. We restrict ourselves to the regular subdivision of a hexahedron into eight child hexahedra. The grid adaptation is subject to a 1-level transition constraint between cells which share a common face or a common edge, see Figure 1. Any refinement technique which observes these rules would be appropriate, as for instance AMR [?], the cartesian refinement based on saturated error indicators [?, ?] or multi-scale analysis [?].

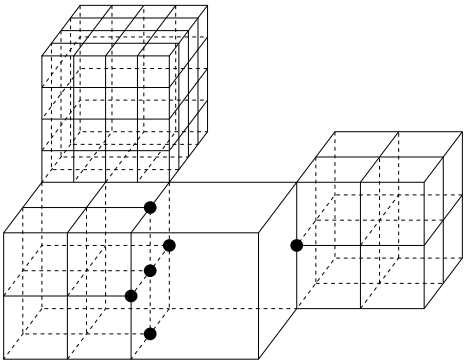

Figure 1: 1-level transition constraint for the primal 3D-grid

Given a primal grid $G_p$ we are now looking for a corresponding staggered *dual grid* $G_d$. The numerical algorithms, and in particular those of finite volume type, require at least:

- $G_p$ and $G_d$ are meshes for the same computational domain $\Omega$.
- Interior nodes, edges and faces of $G_p$ and $G_d$ do not coincide.

Moreover, the dual grid $G_d$ should locally reflect the resolution of the primal grid $G_p$. With respect to the limitation of the timestep size by the CFL-number and the distance of discontinuities in the numerical solution one should also try to maximize the distance between faces of primal and dual cells.

Before we present our new approach for the construction of the staggered grid, we first discuss the well known "diamond grids" (cf. [?, ?]). For reasons of clarity only the 2-dimensional figures are presented. All arguments are formulated also for 3-dimensional grids. The steps to construct dual cells read as follows:

1. Construct edges of dual cells simply by connecting the midpoint of a primal cell $C$ with all of $C$'s corners (see Figure 2(a)). In 2D these edges subdivide $C$ into four triangular dual cell parts, whereas in 3D the dual edges span 12 dual faces which bound six pyramidal dual cell parts.

2. Dual cells result by sticking together adjacent dual cell parts with a common primal face (see Figure 2(b)).

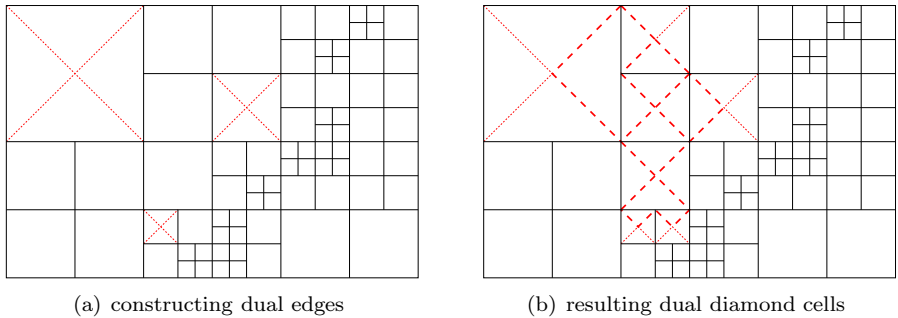

(a) constructing dual edges    (b) resulting dual diamond cells

Figure 2: diamond-type dual cell construction in 2D

This construction technique is impressively simple, but suffers from excessive numerical cost:

- Since every primal grid cell is subdivided into four (in 2D) resp. six (in 3D) dual cell parts, and a dual cell is normally composed by two cell parts the dual grid contains approximately twice resp. thrice as much cells as the corresponding primal grid.

- Following the algorithm of Nessyahu and Tadmor from [?] two kinds of integration have to be performed in (4) to determine the numerical solution on the next time level. The far more expensive part is $I_2$, the integration of flux functions over the boundary of cells. Since the cellwise representation of the numerical data is in general discontinuous at the cell boundaries, the fluxes have to be evaluated in the interior of primal cells (e.g. in the center of gravity of the dual faces in the case of piecewise linear data. Its evaluation in the nodes of the dual cells is not possible since they often coincide with the boundary of primal cells). Since there are four (in 2D) resp. twelve (in 3D) dual faces per primal cell the number of necessary flux evaluations for a timestep onto the dual grid sums up to:

$$\# \text{ fluxes(2D)} = 4|\mathcal{C}_p|$$
$$\# \text{ fluxes(3D)} = 12|\mathcal{C}_p|$$

- Due to the finer resolution of the dual grid the diameter of a dual cell is by factor $\sqrt{2}$ smaller than the diameter of its corresponding primal cell. This leads to a restriction of the global timestep size in the explicit time integration procedure.

Our less expensive alternative approach consists of the following principles:

1. Construct dual cells always surrounding exactly one node $N_p \in \mathcal{N}(G_p)$. Hence dual nodes Since in a cartesian grid there are approximately as many nodes as cells, primal and dual grid are of comparable size :

$$\# \text{ fluxes(2D)} = 2|\mathcal{N}_d| \approx 2|\mathcal{C}_p|$$
$$\# \text{ fluxes(3D)} = 3|\mathcal{N}_d| \approx 3|\mathcal{C}_p|$$

   Moreover, the dual grid reflects to local resolution of the primal grid.

2. The shape of the dual cells is determined by a Voronoi decomposition of the domain $\Omega$ respecting all nodes $N_p \in \mathcal{N}(G_p)$. Nodes of dual cells now lie in the interior of primal cells. For the flux-integration the flux-functions can be evaluated there with respect to every spatial direction. Thus the numerical effort is essentially determined by the number of nodes in the dual grid (which is only slightly larger than the number of primal cells).

The cartesian structure of the primal grid and its adaptation constraints (regular refinement, 1-level transition) lead to a rigidly prescribed location of the primal nodes. This allows the Voronoi decomposition to be performed *locally* on every cell $C_p$ of the primal grid (as it will be outlined in the next chapter). Depending on the set of primal nodes on the boundary of a primal cell $C_p$, $\mathcal{N}^+(C_p) := \mathcal{N}(G_p) \cap C_p$, we deduce $C_p$'s decomposition into local Voronoi regions. All these local regions together form a *local pattern* on $C_p$. Local patterns on adjacent cells of the primal grid automatically match. All local Voronoi regions corresponding to a fixed node $N_p$ of the primal grid finally form the dual cell $C_d(N_p)$.

Due to the primal grid structure and the 1-level transition constraint the number of *essentially different* local patterns is finite, and every possible decomposition can be analysed in advance. Later, the numerical scheme gets instant access to these predefined patterns and uses them in scaled and rotated copies. These copies are retrieved at run time and do not have to be stored.

## Voronoi regions: choice of the norm

Given a set $\Omega \subset \mathbb{R}^n$ and a finite set of nodes $\mathcal{N} = \{N_k, k = 1, \ldots, m \,|\, N_k \in \Omega\}$, we decompose $\Omega$ into $m$ Voronoi regions $V_{N_k}$ by the following conditions:

(i) $\bigcup_{k=1}^{m} V_{N_k} = \Omega$

(ii) $\dot{V}_{N_i} \cap \dot{V}_{N_j} = \emptyset, i \neq j$

(iii) $\|x - N_k\| \leq \|x - N_j\|, x \in V_{N_k}, \forall j \neq k$.

Here $\dot{V}_{N_i}$ denotes the interior of $V_{N_i}$, and $\|\cdot\|$ an arbitrary norm on $\mathbb{R}^n$.

Let us now discuss the choice of an appropriate norm in 2D more in detail. The same arguments hold also for 3D.

A Voronoi decomposition of $\Omega \subset \mathbb{R}^2$ refering to *only two* nodes $N_i$ resp. $N_j \in \Omega$ splits $\Omega$ into two Voronoi regions, $V_{N_i(N_j)}$ (around $N_i$) and $V_{N_j(N_i)}$ (around $N_j$), sharing a common *separating polygon* $S_{ij}$. This separating polygon is the locus of the intersection points of circles with same diameter centered at $N_i$ resp. $N_j$. The shape of the separating polygon depends both on the locus of $N_i$ and $N_j$ as well as on the used norm on $\mathbb{R}^2$. In general, a Voronoi region $V_{N_i}$ around the node $N_i$ takes the form

$$(5) \qquad V_{N_i} = \bigcap_{N_j \in \mathcal{N}} V_{N_i(N_j)}$$

Respecting the ordinary $\|\cdot\|_2$-norm, any $S_{ij}$ is a straight line, hence Voronoi regions $V_{N_i}$ are always convex. In contrast, for the $\|\cdot\|_\infty$-norm, $S_{ij}$ consists, in general, of three straight lines, which are aligned with the cartesian axes ($x$-axis, $y$-axis) or the plane diagonals (of the $xy$-plane), see Figure 3(a). For special loci of $N_i$ and $N_j$ the separating polygon simplifies to a straight line (see Figure 3(b)), or becomes not-unique in regions away of $N_i$ and $N_j$ (see Figure 3(c)). For the latter case we choose the prolonged straight line from the unique part as the separating polygon. Moreover, Voronoi regions $V_{N_i}$ might happen to be non-convex.

However, thanks to the cartesian structure of the primal grid the $\|\cdot\|_\infty$-norm even simplifies the construction of local Voronoi regions on cells of the primal grid. Compared to the corresponding construction for the Euclidean norm, boundary faces of Voronoi regions are now aligned only with the axes or the plane diagonals of the primal grid (cf. Figure 4(b) vs. Figure 4(a)). We therefore favour the $\|\cdot\|_\infty$-norm.

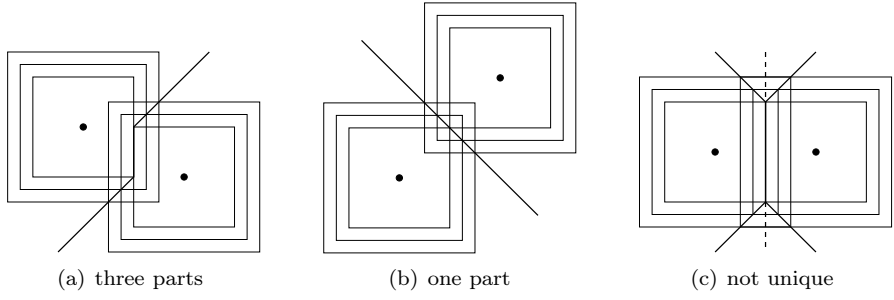

(a) three parts  (b) one part  (c) not unique

Figure 3: Different configurations of two points, some $\|\cdot\|_\infty$-norm circles around these points, and the resulting separating polygons in 2D

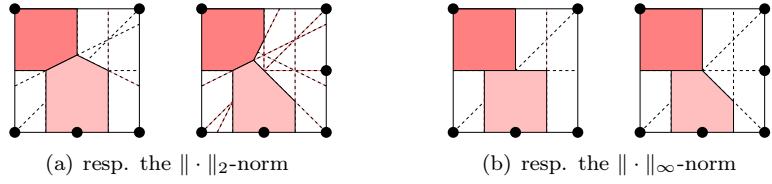

(a) resp. the $\|\cdot\|_2$-norm  (b) resp. the $\|\cdot\|_\infty$-norm

Figure 4: Comparison of Voronoi regions and separating polygons on a primal cell in 2D

## Local Voronoi construction

Next, we inverstigate the locality of the Voronoi construction. With (5) we get bounding boxes $V_{N_i}^{C_p}$ of local Voronoi regions $V_{N_i}$ on a primal grid cell $C_p$ by

$$(6) \qquad V_{N_i}^{C_p} := C_p \cap \bigcap_{N_j \in \mathcal{N}(C_p)} V_{N_i(N_j)} \supseteq C_p \cap V_{N_i}$$

$\mathcal{N}(C_p)$ denotes the set of corners of the cell $C_p$. Figure 5 illustrates the bounding boxes resulting by intersection (6) on a primal cell in 3D. Here, $N_i$ is a corner node (fig. 5(a)), a node in the midpoint of an edge (fig. 5(b)) resp. a node in the midpoint of a face (fig. 5(c)) of $C_p$.

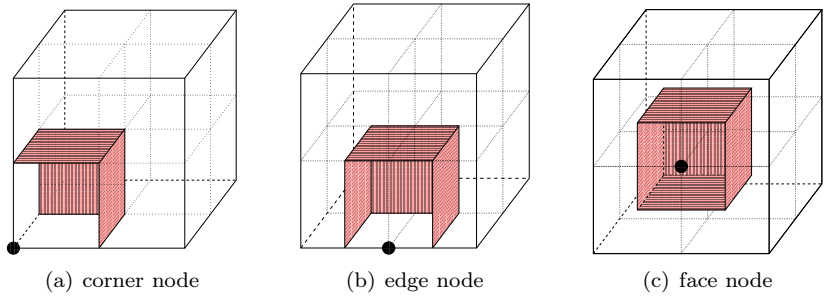

(a) corner node  (b) edge node  (c) face node

Figure 5: Local Voronoi construction in 3D

In order to get the final shape of a local Voronoi region $V_{N_i}$ on $C_p$ it suffices to perform the intersection (5) for all nodes $N_j \in \mathcal{N}^+(C_p)$, i.e. respecting only the nodes on the boundary of $C_p$. Shorter, we state

**Lemma 1.**
$$C_p \cap V_{N_i} = C_p \cap \bigcap_{N_j \in \mathcal{N}^+(C_p)} V_{N_i(N_j)}$$

*Proof.* Assuming lemma 1 were false. Then there is a node $N_k \in \mathcal{N} \setminus \mathcal{N}^+(C_p)$ such that $V_{N_i}^{C_p} \cap V_{N_i(N_k)} \subsetneq V_{N_i}^{C_p}$ which implies

(7) $$\exists\, x \in V_{N_i}^{C_p} \text{ with } \|x - N_k\|_\infty < \|x - N_i\|_\infty$$

i.e. $x$ should not be assigned $V_{N_i}$.

Let $\mathfrak{H}_{N_i}(C_p) := \bigcup_{y \in V_{N_i}^{C_p}} B_{\|y-N_i\|_\infty}(y)$. Now, (7) implies

(8) $$N_k \in \mathfrak{H}_{N_i}(C_p)$$

With $\Delta x = \text{diam}(C_p, \|\cdot\|_\infty)$ we define $\mathfrak{U}_{N_i}(C_p) := \bigcup_{y \in V_{N_i}^{C_p}} B_{\frac{\Delta x}{2}}(y)$ and

$$\mathfrak{H}_{N_i}^+(C_p) = \begin{cases} \mathfrak{U}_{N_i}(C_p) \setminus \hat{C}_p(N_i) & \text{if } N_i \text{ is corner node of } C_p \\ \mathfrak{U}_{N_i}(C_p) & \text{else} \end{cases}$$

where $\hat{C}_p(N_i)$ denotes that neighbour of $C_p$ which shares only the node $N_i$, i.e. $C_p \cap \hat{C}_p(N_i) = N_i$.

Suppose that we are in case 5(a). Let $z \in \hat{C}_p(N_i)$ be fixed. Then we have $\|y - N_i\|_\infty < \|y - z\|_\infty$ for any $y \in V_{N_i}^{C_p}$. Hence there is no $y \in V_{N_i}^{C_p}$ such that $z \in B_{\|y-N_i\|}(y)$, thus $\hat{C}_p(N_i) \cap \mathfrak{H}_{N_i}(C_p) = \emptyset$.

In general, since $\|y - N_i\|_\infty \leq \frac{\Delta x}{2}\ \forall y \in V_{N_i}^{C_p}$, the relation $\mathfrak{H}_{N_i}(C_p) \subseteq \mathfrak{H}_{N_i}^+(C_p)$ holds for every case 5(a), 5(b) and 5(c). Due to the 1-level grading condition over faces and edges of $C_p$ the only grid nodes in the interior of $\mathfrak{H}_{N_i}^+(C_p)$ are the possibly occuring hanging nodes on $C_p$'s boundary. Hence

(9) $$\mathfrak{H}_{N_i}^+(C_p) \cap \mathcal{N} \setminus \mathcal{N}^+(C_p) = \emptyset$$

i.e. there is no $N_k$ fulfilling condition (8). This proves lemma 1. □

**Remark 1.** *The special treatment of the case 5(a), where $N_i$ is a corner of $C_p$, is necessary since the 1-level grading condition over faces and edges still allows a 2-level transition over nodes. In this case $\mathfrak{U}_{N_i}(C_p)$ and $\mathcal{N} \setminus \mathcal{N}^+(C_p)$ would share a node of $\hat{C}_p(N_i)$, and (9) would not hold.*

**Remark 2.** *The arguements in the proof do not depend on the spatial dimension of $C_p$. For a better understanding the reader might prefer to draw sketches of the proof in the 2-dimensional setting.*

## 3   Pattern construction in 2D

In contrast to the local patterns from Figure 9(a) used for the dual grid construction in [?] we now want to create them by a clear geometrical rule. By means of a local Voronoi decomposition respecting the maximum norm we assign the "volume" of a cell $C_p$ of the primal grid the nodes $N_k \in \mathcal{N}^+(C_p)$ on the boundary of $C_p$. These are the four corner nodes of $C_p$ (they always exist) and possibly up to four nodes on the midpoint of $C_p$'s edges. These additional hanging nodes only exist if adjacent cells in the primal grid were refined. The local intersection of seperating polygons boils down to six different types of local Voronoi regions, cf. Figure 6.

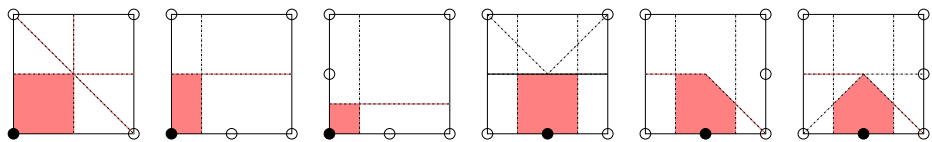

Figure 6: Resulting local Voronoi regions on a primal cell in 2D

Obviously, all cuts of the cell $C_p$ run parallel to the axes or the plane diagonals. We could therefore describe the steps for the construction of local Voronoi regions on $C_p$ equivalently by the following algorithm:

1. Subdivide the cell $C_p$ into 16 congruent *squares*.

2. Subdivide all squares containing the midpoint of $C_p$ into two *triangles*. The cuts follow the face diagonals through $C_p$'s midpoint. On $C_p$ we get 8 triangles.

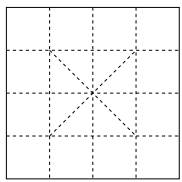

Figure 7: Subdivision of a primal cell in 2D

   These squares and triangles, depicted in Figure 7, will be denoted as *atoms* of $C_p$.

3. Finally, we assign these 20 atoms of $C_p$ the nodes of $\mathcal{N}^+(C_p)$ simply by calculating the $\|\cdot\|_\infty$-distance between the nodes and the atom's center of gravity. All atoms being assigned the same node $N_k \in \mathcal{N}^+(C_p)$ form the *local Voronoi region* $V_{N_k}(C_p)$ on $C_p$. The set of all local Voronoi regions on $C_p$ forms the *local pattern*. Figure 8 shows all six essentially different local patterns in 2D.

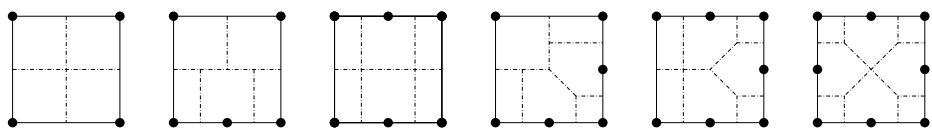

Figure 8: All essentially different local patterns in 2D

Assembling scaled and rotated copies of appropriate local patterns lead to the dual grid depicted in Figure 9(b). Compared to the dual grid from Figure 9(a) proposed in [?] there is now a unique correspondance between dual cells and primal nodes. This simplifies the adressing of data on the dual grid.

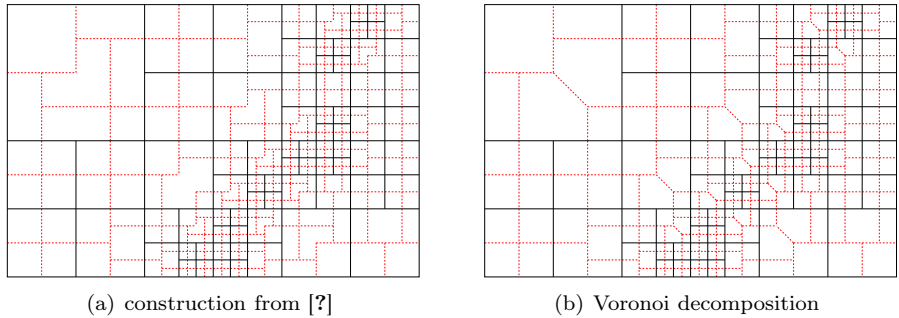

(a) construction from [?]         (b) Voronoi decomposition

Figure 9: 2D cartesian primal grid (solid) and corresponding dual grid (dashed)

## 4   Pattern construction in 3D

Now we extend the steps to create local patterns described in the previous section to 3-dimensional grids. Applying the concept of constructing local Voronoi regions in 3D leads to an only slightly more complicated description than in 2D. Again, we assign the volume of a primal cell $C_p$ the nodes on its boundary, i.e. all $N_k \in \mathcal{N}^+(C_p)$. These are the eight corner nodes of $C_p$, and possibly up to 18 nodes on the midpoint of $C_p$'s six faces and 12 edges. These hanging nodes occur if an adjacent cell, sharing a common face resp. edge with $C_p$, was refined. For a general location of two nodes $N_i$ and $N_j$ in $\mathbb{R}^3$ the *separating surface* between the Voronoi regions $V_{N_i(N_j)}$ and $V_{N_j(N_i)}$ consists of up to seven planar patches, whose normals are aligned with the axes ($x$-, $y$-, $z$-axis) or the plane diagonals ($xy$-, $yz$-, $zx$-plane).

Fortunately, the location of nodes in the primal grid restricts the variety of occuring separating surfaces to the six cases of Figure 10. Finally, on a primal cell only scaled, translated and rotated versions of these cuts have to be executed to get the local Voronoi regions.

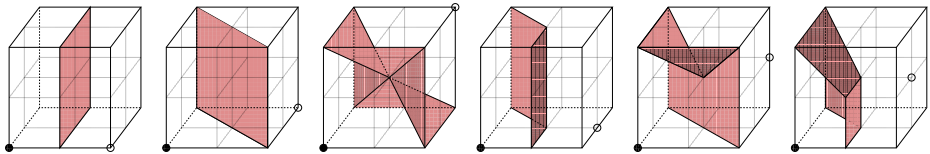

Figure 10: Occuring separating surfaces in 3D

Hence the whole construction can again be described as a subdivision algorithm:

1. Subdivide the hexahedral cell $C_p$ into 64 congruent *cubes*.

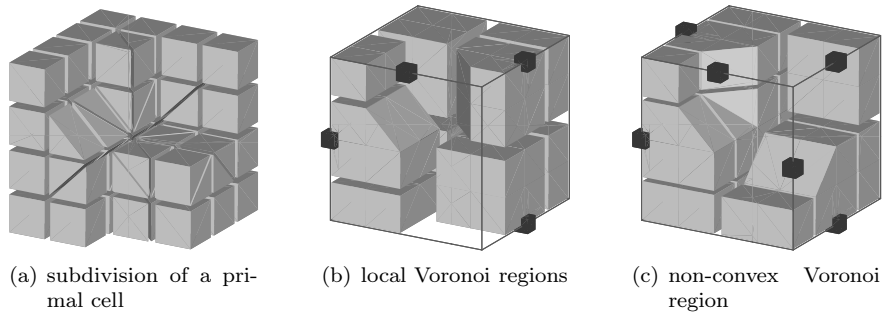

(a) subdivision of a primal cell    (b) local Voronoi regions    (c) non-convex Voronoi region

Figure 11: Construction of local Voronoi regions in 3D

2. Subdivide all cubes containing a midpoint of $C_p$'s faces into two *prisms*. The cuts follow the face diagonals through the face midpoints. On $C_p$ we get $8 \times 6 = 48$ prisms.

3. Subdivide all cubes containing the central point of $C_p$ into six *tetrahedra*. The corresponding three cuts equal the diagonal cuts on their adjacent cubes from step 2. This results in $6 \times 8 = 48$ tetrahedra.

These cubes, prisms and tetrahedra will again be denoted as *atoms* of $C_p$. Figure 11(a) shows a cell $C_p$ and some of its atoms. Next we assign these $64+48-24+48-8 = 128$ atoms of $C_p$ the nodes of $\mathcal{N}^+(C_p)$ by calculating the $\|\cdot\|_\infty$-distance between the nodes and the atom's center of gravity. All atoms being assigned the same node $N_k \in \mathcal{N}^+(C_p)$ form the *local Voronoi region* $V_{N_k}(C_p)$ on $C_p$. Figure 11(b) shows a cell $C_p$, the set $\mathcal{N}^+(C_p)\setminus\mathcal{N}(C_p)$ (i.e. the corners of $C_p$ left out) and a few corresponding local Voronoi regions. Let us mention that local Voronoi regions might happen to be slightly non-convex (region in the center of Figure 11(c)). But this should not effect the feasibility of our approach, and we do not expect additional problems in constructing a second order central finite volume scheme by the occurrence of these dual cells.

The set of all local Voronoi regions on $C_p$ form the *local pattern*. The shape of these three-dimensional local patterns on the boundary faces of a primal cell matches the two-dimensional patterns from Figure 8.

## Counting local patterns

Since the number of different distributions of nodes on the boundary of a primal cell is finite, all local patterns can be assembled and stored in advance. Later, running the numerical scheme, these patterns are used in scaled and rotated copies. To determine the number of essentially different local patterns (that do not arise from each other by rotation or reflection) we used tools of the Pólya theory (see [?], [?]). All information for the counting problem is provided by the group of symmetries $\mathcal{D}(\mathfrak{W})$ of the cube $\mathfrak{W}$. Analysing how the elements of $\mathcal{D}(\mathfrak{W})$ act on the set of $\mathfrak{W}$'s edges and faces leads to the number of 227 essentially different local patterns. For details we refer to [?]. The list of all 227 patterns (some of them with figures) can be found on http://www.igpm.rwth-aachen.de/wolfram/local_patterns.html

They can all be assembled, analyzed and stored in a lookup-table. Every local refinement situation can be identified with exactly one *reference pattern*. This reference pattern as well as an appropriate mapping (i.e. permutation matrix) onto the really occuring pattern are a priori known.

# 5   Applications

At the current stage, our implementation of the finite volume scheme includes the primal and dual grid generation, the integration of piecewise linear functions over volumes and boundary faces of cells, but still lacks some essential components (like data reconstruction and limiting, treatment of divers boundary conditions) to end up with a second order accurate non-oscillatory method. Nevertheless, we have verified our dual-grid approach with standard test cases.

After summarizing the algorithmical demands and presenting our integration strategy, we address the data access, describe our test cases and finally compare the expected numerical effort of our method with the diamond grid approach and the established non-staggered HLL-scheme.

### Algorithmical demands

The evaluation steps in (4) are twofold. Part $I_1$ calls for an integration of the piecewise linear function $v^n$ over a bounded volume and can be evaluated exactly. Part $I_2$ is a bit more intricate. For an illustration we refer to Figure 12. The normal part of a non-linear flux-function $F$ has to be integrated in time over an interval $\Delta t = t^{n+1} - t^n$ and in space over a polygonally bounded planar cell face $\mathfrak{p}$. In general, a face $\mathfrak{p}$ of a cell $C^*$ consists of several polygonal parts $\mathfrak{p}_i$ in adjacent cells $C_i$ of the corresponding staggered grid.

Any second order accurate quadrature formula for $I_2$ should evaluate $F$ in time and space where its argument $v = v(x,t)$ is uniquely defined. For the temporal part of $I_2$ we use the midpoint rule. For that reason $v$ has to be extrapolated to approximate $v(x, t^{n+1/2})$. As time evolves, characteristics spread out from the boundary into the interior of each cell. In order to avoid the solution of local Riemann problems at the cell interfaces the spatial integration in $I_2$ should choose quadrature points away from the boundaries which are, at time $t^{n+1/2}$, not yet reached by the characteristics. The distance of these points to the cell boundary limit the size of the global timestep $\Delta t = t^{n+1} - t^n$.

Though admissible, the choice of quadrature points in the center of gravity of every polygonal part $\mathfrak{p}_i$ is disadvantageous:

- In a cartesian 3D-grid there about three times more faces than cells. A polygonal face contains of mostly four polygonal parts. We would have to evaluate one flux function on every polygonal parts.

- The distance of the center of gravity of a polygonal part $\mathfrak{p}_i$ to the boundary of the cells $C_i$ of the corresponding staggered grid is rather short. This limits the global timestep size unnecessarily.

In our algorithm we implemented a generalized trapezoidal rule on polygons. The flux function $F$ is evaluated at the nodes of a polygonal face. Later, these values are scaled with barycentric weights respecting the face geometry to achieve a second order accurate integration rule. This approach promises to be more competitive:

- In a cartesian 3D-grid the number of nodes and the number of cells is similar. We evaluate all three directional flux functions at every node $\mathfrak{v}_i$. The numerical effort is by factor $\approx 4$ less than the method described above.

- Since the distance of nodes of a polygonal face $\mathfrak{p}$ of a cell $C^*$ to the boundary of cells $C_i$ from the corresponding staggered grid is larger than the distance from the $\mathfrak{p}_i$'s center of gravity to the boundary the global timestep size is less limited compared to the method above.

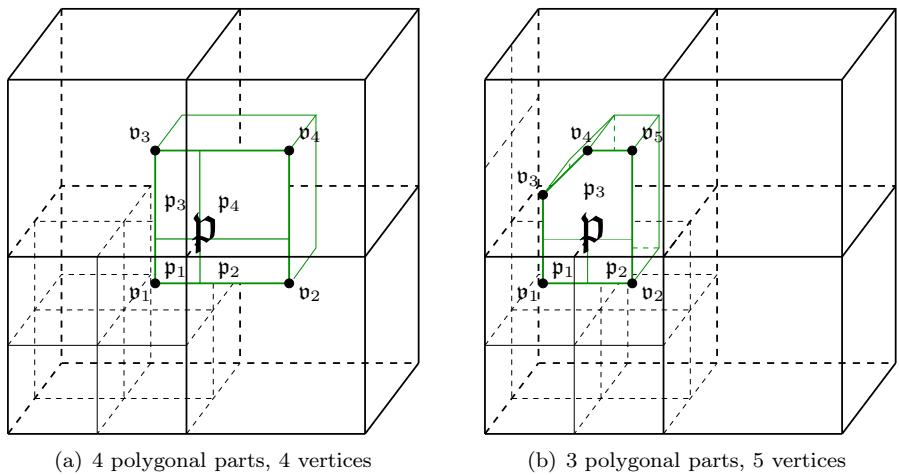

(a) 4 polygonal parts, 4 vertices    (b) 3 polygonal parts, 5 vertices

Figure 12: Compound boundary faces of dual cells

## Data access

For fast access all geometric and numerical data, both for the primal and the dual grid, are stored in (some of them level-linked) hashmaps. On every cell $C$ of the grid $G$ the finite volume scheme needs to know

1. the volume of $C$
2. its neighbouring cells on the grid $G$ (for a higher order reconstruction)
3. the common faces of $C$ and its neighbours, including their normals and their area, as well as the nodes on the grid $G$ that form these faces (to determine fluxes over faces by an appropriate quadrature rule)

On the cells of the cartesian primal grid $G_p$ such information is quite easily accessible by index shift operations. In contrast, for a cell $C_d$ of the dual grid $G_d$ this knowledge is scattered over all local patterns contributing to the construction of $C_d$ and has to be "collected" during a complete traversal of the primal grid $G_p$. Hence finding the appropriate local pattern on a primal cell is an often used operation and has to be fast.

The local pattern on a primal grid cell $G_p$ is determined by the occurrence of hanging nodes in the midpoints of $G_p$'s 12 edges and six faces (cf. Section 4. These information form an 18-bit key which allows instant access to the corresponding reference pattern (one of 227, cf. Section 4) as well as the appropriate permutation matrix from a predefined hashmap. This hashmap can easily be initialized in a pre-roll step by applying each of the 48 self-mappings of the cube to all of the 227 essentially different local refinement constellations. Finally, one gets 6210 hashmap entries.

A local pattern consists of a certain number of local Voronoi regions $V_{N_k}$, each of them composed of several atoms (cubes, prisms and tetrahedra, see Figure 11(a)). We get all information listed above locally:

1. the volume of $V_{N_k}$ equals the sum of the volumes of its atoms

2. neighbouring Voronoi regions on a pattern refer to neighbouring cells on the dual grid $G_d$

3. common faces of neighbouring Voronoi regions are the union of common faces of their atoms. Thus, the normals, areas and the face-forming nodes are known by construction of the local pattern.

Whereas on the primal grid $G_p$ all cells $C_p$ are cubes and common faces to neighbouring cells are always square-shaped, there is a wider range of different cell- and face-types on the dual grid $G_d$, see Figure 13(a) for a rough impression.

## Example I: Volume integration

In order to test the presented dual grid construction and to roughly estimate the numerical cost of a staggered grid scheme, a *volume integration* over a subset of the dual grid $G_d$ was performed. The primal grid has been constructed by adaptively resolving the isosurface of a tilted elliptical paraboloid inside the unit cube, see Figure 13(b).

We verified our algorithm by means of the Gauss integral theorem

$$\int_{\hat{\Omega}} \text{div } \mathfrak{v} \ dx = -\int_{\partial\hat{\Omega}} \mathfrak{v} \cdot \mathfrak{n} \ dx$$

by choosing $\mathfrak{v}(\vec{x}) = \frac{1}{3}\vec{x}$, $\forall \vec{x} \in \mathbb{R}^3$, a subdomain $\hat{\Omega} \subseteq \Omega$, and integrating over all cells $C_d \subset \hat{\Omega}$ of the dual grid $G_d$

1. by summing up the volumes of the involved Voronoi regions

2. by integrating $\frac{1}{3}\langle\vec{x}, \mathfrak{n}_f\rangle$ over all common faces $f$ of neighbouring Voronoi regions (in view of flux integration over cell boundaries). Here $\mathfrak{n}_f$ denotes the oriented normal of $f$.

We performed the volume integration described above for several resolutions of the primal grid $G_p$ and examined

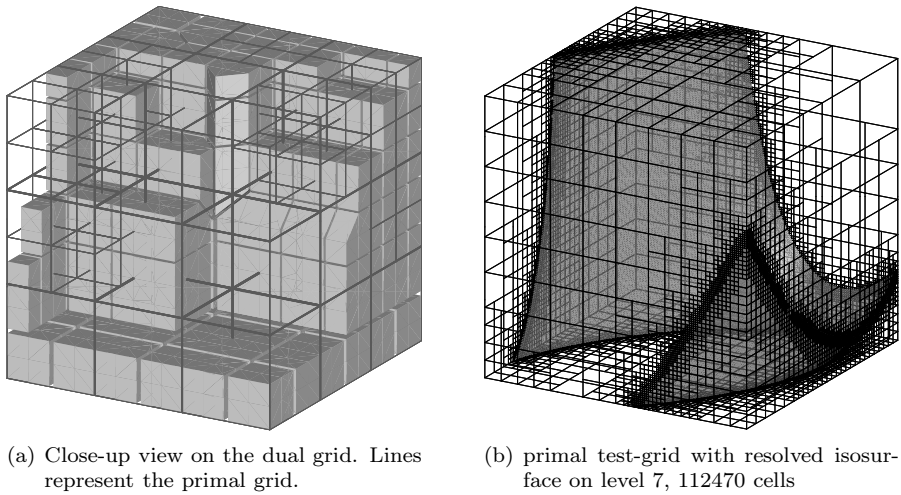

(a) Close-up view on the dual grid. Lines represent the primal grid.

(b) primal test-grid with resolved isosurface on level 7, 112470 cells

Figure 13: Primal and dual grid

- the number of different occuring local patterns:
  it leveled off at about 40

- the frequency of occurrence of these patterns:
  approximately 80% of the cells $C_p$ on the primal grid refer to the most simple pattern (no additional node on the boundary of $C_p$), whose numerical cost are low

- the expected number of flux evaluations for the staggered grid solver by counting the number of nodes in the dual grid

To estimate the profit of our approach, the numerical effort was compared with a cheap non-staggered finite volume scheme where flux integrals are approximated by the Harten-Lax-van Leer (HLL) Riemann solver (cf. [?]). Results are listed in table 1 below. The numerical cost of both schemes are similar.

| level | primal cells | dual cells | fluxes non-staggered | fluxes staggered |
|---|---|---|---|---|
| 4 | 568 | 1019 | 3810 | 4779 |
| 5 | 4502 | 6969 | 30126 | 33717 |
| 6 | 24781 | 34710 | 163899 | 161001 |
| 7 | 112470 | 149011 | 737283 | 690420 |
| 8 | 493368 | 630364 | 3210945 | 2906742 |
| 9 | 2156547 | 2684250 | 13944345 | 12220641 |

Table 1: Comparison between staggered and non-staggered approach

## Example II: Rotating Cone

Not yet having incorporated linear data reconstruction and limiting in our method (what would be essential for a second order non-oscillatory scheme) we run a standard linear advection problem, the "Rotating cone". The comparison to the diamond-grid approach as well as to the non-staggered HLL-scheme in terms of numerical effort confirms our topological statements from Section 2 and encourages for the further work.

**Initial data & numerical solution**

As initial data for the rotating cone advection problem we resolve a smooth function $f$ with compact support over an octant of the surface $\partial B$ of a ball $B = B(c, R)$ (cf. Figure 14(a)) The ball $B$ is centered at $c = (c_x, c_y, c_z) = (0.6, 0.3, 0.2)$ and has a radius $R = 1/4$. The function $f$ is defined as follows:

$$\begin{aligned}
r^2 &= (x - c_x)^2 + (y - c_y)^2 + (z - c_z)^2 \\
q &= 4|1 - r^2| \\
f(q) &= \begin{cases} 1 - 2q^2 & \text{for } q < 1/2, x, y, z \geq 0 \\ 2(q - 1)^2 & \text{for } 1/2 \leq q \leq 1, x, y, z \geq 0 \\ 0 & \text{else} \end{cases}
\end{aligned}$$

This cone is now rotated around $c$ in the plane spanned by $v = (1, 0, 0)^t$ and $w = (0, 1, 0.5)^t$. The advection performs up to time $t = \pi/4$ at constant angle velocity 1. The exact solution is depicted in Figure 14(b). Since our scheme is still only first order accurate the numerical solution suffers from immoderate smearing (cf. Figure 14(c)).

**Grid specifications & numerical effort**

Nevertheless, the resulting grids offer an opportunity to estimate and to compare the numerical effort of our scheme. The method to count the number of necessary flux evaluations in two successive timesteps is outlined in table 2 and performed for the grids from Figure 14(c) in table 3.

| timestep | Diamond | HLL | Voronoi |
|---|---|---|---|
| primal → dual | 12 #primal cells | 2 #primal faces | 3 #dual nodes |
| dual → primal | 1 #primal faces | 2 #primal faces | 3 #primal nodes |

Table 2: Counting flux evaluations in two successive timesteps

| Grid specifications | |
|---|---|
| primal cells | 112106 |
| primal faces | 345564 |
| primal nodes | 121561 |
| dual nodes | 152630 |

| # Flux evaluations | | | |
|---|---|---|---|
| | $\mathbf{G}_p \to \mathbf{G}_d$ | $\mathbf{G}_d \to \mathbf{G}_p$ | $\sum \#f$ |
| Diamond | 1.35 M | 0.35 M | 1.69 M |
| HLL | 0.69 M | 0.69 M | 1.38 M |
| Voronoi | 0.46 M | 0.36 M | 0.82 M |

Table 3: Grid specifications and numerical effort

Obviously, in an adaptively octtree-refined cartesian 3D-grid the number of cells on the (locally) deepest refinement level, i.e. those without hanging nodes (we call it *property* $\mathcal{H}$), outranges the number of coarser cells by far. In general, discontinuities in the numerical solution are resolved by (possibly thin) 2-dimensional manifolds. For a rough argument, let $N_L$ be the number of grid cells and $A_L$ the number of cells with property $\mathcal{H}$ on a grid with at most $L$ refinement levels. Resolving a discontinuity by refining $k$ cells on refinement level $L$ effects the numbers $A$ and $N$ as follows: each of the $k$ cells generates 8 child cells, the original cell itself disappears, thus $N_{L+1} \approx N_L + 7k$. The refinement of a cell inserts hanging nodes on its coarser neighbour cells. Since all $k$ cells are supposed to form a thin surface, each of them "destroys" the property $\mathcal{H}$ only in the two face-neighbours in the surface's normal direction. The refined cell itself contributed most likely to $A_L$. Every child cell on level $L+1$ obviously has property $\mathcal{H}$. Thus $A_{L+1} \approx A_L + 5k$. For $k \gg 1$ one gets the ratio $r = A/N \approx 5/7$. If the $k$ cells on the fine grid level rather form a volume than a thin surface, the ratio $r$ gets even closer to 1.

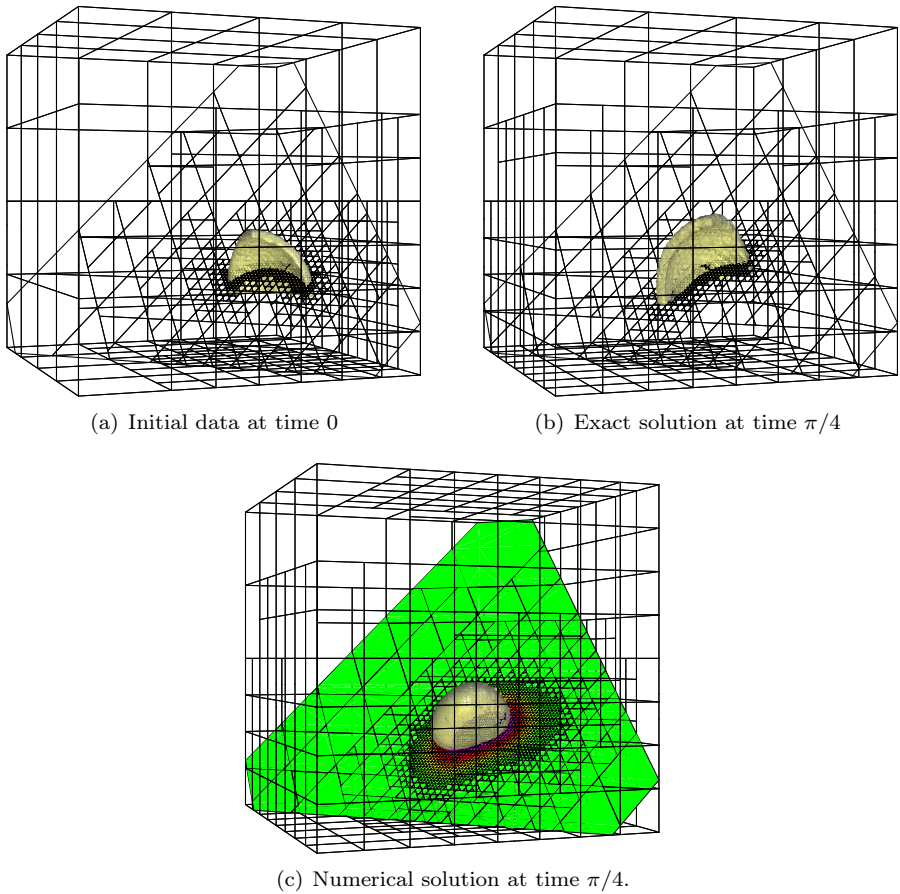

(a) Initial data at time 0

(b) Exact solution at time $\pi/4$

(c) Numerical solution at time $\pi/4$.

Figure 14: Rotating cone with first order scheme

In this regard, the ratio of 94% in our example from Figure 14(c) does not surprise. For all these cells the according local pattern is the easiest possible, only containing one dual node which leads to only half as much local costs than incurred by the compared HLL-solver, and only one fourth in comparison to the diamond cell approach.

## 6 Conclusions & extensions

In the current paper basic steps of the implementation of Nessyahu and Tadmor's staggered scheme [?] on adaptively refined cartesian grids in 3D were introduced. The construction of dual grid cells as Voronoi regions respecting the maximum-norm simplifies their geometrical description and, in the end, the implementation. Substantial combinatorial investigation in advance allows our scheme in return to exploit local geometrical structure of the grid at run-time to save a lot of numerical effort.

Though up to now only first order, the extension of the staggered grid approach to a higher order 3D finite volume scheme is in process. It concerns the piecewise linear reconstruction and limitation of cell- and flux-values to end up with a second order scheme, and the data handling for large 3D simulations. Still challenging are the incorporation of obstacles in the computational domain and the treatment of divers boundary conditions, especially when general CAD-geometries (and not only cartesian grid cells) need to be incorporated. Additionally desirable would be the parallelization of the staggered grid scheme as well as the local adaptivity in time.